\documentclass[a4paper,10pt,oneside]{article}

\usepackage{tgpagella}
\usepackage[T1]{fontenc}
\usepackage{latexsym}
\usepackage{amssymb}
\usepackage{amsmath}

\usepackage{theorem}

\usepackage{hyperref}

\usepackage{pgf}
\usepackage{xxcolor} 

\usepackage{natbib}

\theoremheaderfont{\scshape}
\theorembodyfont{\rmfamily}

\bibpunct[, ]{(}{)}{,}{a}{}{;}

\newtheorem{theorem}{Theorem}

\def\qed{\relax
   \ifmmode
    ~\hfill\Box
   \else
    \unskip\nobreak ~\hfill$\square$%
   \fi \par}

\newcommand{\sep}{$\cdot$ }

\begin{document}
	\title{
	The supertask of an infinite lottery}
	\author{Timber Kerkvliet}
	\date{}
	\maketitle

\begin{abstract}
We mathematically model the supertask, introduced by \citet{hansen,hansen2}, in which an infinity of gods together select a random natural number by each randomly removing a finite number of balls from an urn, leaving one final ball. We show that this supertask is highly underdetermined, i.e. there are many scenarios consistent with the supertask. In particular we show that the notion of uniformity for finitely additive probability measures on the natural numbers emerging from this supertask is unreasonably weak.
	
\end{abstract}

\section{Introduction}
\label{sec:introduction}

A supertask is a task that consists of countably many tasks that are carried out within a finite time interval. Even if every task of a supertask is unambiguously described, it is possible that it is not uniquely determined what will happen when the supertask is performed. In general, to determine which scenarios are consistent with a description of a certain supertask, we need a (mathematical) model of the supertask that specifies all the (a priori) possible scenarios and the constraints on those scenarios as dictated by the supertask. Within such a model, one can try to determine the set of all scenarios consistent with the supertask, i.e. the set of possible scenarios that satisfy the contraints.

\citet{hansen,hansen2} introduces a supertask in which an infinity of gods together select a random natural number by each randomly removing a finite number of balls from an urn, leaving one final ball. Every task consists of randomly removing a finite number of balls and hence is probabilistic in nature. This means that by the very nature of the supertask, what will happen is undetermined. However, the probability of what happens in every task, given what has happened up to that task, is completely determined. Therefore, we consider finitely additive probability measures of outcomes, which we call \emph{probability functions} from now on, as scenarios rather than outcomes themselves and then ask the question what scenarios are consistent with the supertask. (In Section \ref{sec:story} we explain why we use probability functions instead of probability measures)

At face value, it seems perfectly reasonable to use this supertask to give a definition of `uniformity' for probability functions on the natural numbers, by calling all probability functions of the final ball consistent with the supertask `uniform'. Hansen conflates this notion of uniformity that emerges from the supertask itself, with notions of uniformity proposed on completely different grounds by mathematicians.\footnote{Examples of notions of uniformity proposed by mathematicians are found in \citet{definetti}, \citet{schirokauerkadane}, \citet{wenmackershorsten} and \citet{kerkvlietmeester}.} It is \emph{a priori} not clear that any of these proposed notions is equivalent with the notion of uniformity that emerges from the supertask. A mathematical model followed by an analysis of that model, is necessary to conclude that.

In this paper, we present a mathematical model for the supertask of Hansen to study the notion of uniformity that emerges from the supertask. Within this model, we only insist on the finite additivity of candidate probability functions. We prove that given any $p \in [0,1]$ and set $A \subseteq \mathbb{N}$ such that $A$ and $A^c$ are infinite, there is a probability function consistent with the supertask that assigns probability $p$ to the final ball being in $A$. We also show that if $A$ or $A^c$ is finite, then every probability function consistent with the supertask has to assign respectively probability $0$ or $1$ to the final ball being in $A$.

Our analysis shows that the description of the supertask is highly underdetermined: there are many probability functions consistent with the supertask, all giving different probabilities of the final ball being in certain sets. As a consequence, we cannot speak of this supertask as having an uniquely determined underlying probability function. The underdetermination makes the notion of uniformity emerging from the supertask unreasonably weak: it, for example, violates that every residue class modulo $m$ (e.g. in case $m=2$ the odd and even numbers) has probability $1/m$. Hence it casts doubt on calling the supertask a `uniform' lottery on $\mathbb{N}$. 

The structure of this paper is as follows. In Section \ref{sec:story} we give our own short description of the supertask of Hansen and describe in words the results we present in Section \ref{sec:model}, including some intuitions behind the proof. This enables a reader who is not familiar with the mathematical concepts to skip Section \ref{sec:model}, which contains the mathematical result in full technical detail. In Section \ref{sec:discussion} we give our reflection on the results.

\section{The supertask}
\label{sec:story}

We start to give our own compromised description of the supertask:
\begin{quote}
An infinite number of gods, all indentified by an unique natural number, gather to hold a lottery. To do that, for every natural number $k$, they produce a ball of size $2^{-k}$ (so that the total volume of the balls is finite) with the number $k$ and put them in an urn of finite size. Then they follow the following instruction: at time $t_k = 1/k$ god $k$ takes the urn. If god $k$ finds $k+1$ balls in the urn, he randomly removes one. If god $k$ not finds $k+1$ balls in the urn, he first empties the urn and then puts balls $1$ to $k$ in the urn.\footnote{The explicit instruction for the case that god $k$ finds not $k+1$ balls is given to ensure that this part of the instruction is actually never used: since the instruction is also followed by god $k+1$, god $k$ necesarrily finds $k+1$ balls in the urn.}
\end{quote}

\noindent  In the version of Hansen god $k$ removes $2^{k-1}$ balls. We simplified this by letting each god remove only one ball. This only changes the supertask cosmetically, because the only crucial part is that every god removes a \emph{finite} number of balls. More precisely, the task of god $k$ in the description of Hansen (reducing the number of balls from $2^k$ to $2^{k-1}$) is precisely the task that gods $2^{k}-1$ to $2^{k-1}$ perform together in our description. The main reason for this adjustement is that it turns out to simplify the mathematics in Section \ref{sec:model}.

We write an outcome of the process as $\{3\},\{3,7\},\{3,4,7\},\ldots$, where every set in the sequence adds one (new) natural number to its precursor and every natural number is added at some point in the sequence. The first set represents the set of balls after god $1$ has removed a ball, followed by the set of balls after god $2$ has removed a ball, etcetera. So in this case god $2$ found balls $3,4$ and $7$ and removed ball $4$. After that, god $1$ removed ball $7$ leaving ball $3$ as the final ball. On the outcome space consisting of all such sequences, we want to define our probability function. The constraint on this probability function that we get from the supertask is that for every $k$, the conditional probability of the set of balls after god $k$ being in some collection $S$ (consisting of sets of size $k$), given that after god $k+1$ there are precisely $j$ balls that upon removal give an element of $S$, is $j/(k+1)$.

We have two arguments for only requiring a probability function instead of a probability measure,  which is $\sigma$-additive. In the first place, we think a probability function on an infinite space should not necessarily have to assign zero probability to a disjoint union of infinitely many sets of probability zero, if the union is so big that it has the same cardinality as the outcome space. This is, for example, a crucial property of Lebesgue measure on $[0,1]$: every single point has probability zero, but intervals, while being unions infinitely many of points, have positive probabilities.\footnote{\label{foot:lebesgue}Notice that Lebesgue measure on $[0,1]$ does not assign a positive probability to \emph{every} uncountable set: the Cantor set has probability zero while being uncountable. The point here is that it assigns positive probability to \emph{some} uncountable sets. } Since $\sigma$-additivity forces all disjoint unions of  \emph{countably} many sets of probability zero to be zero, we can only guarantee the preservation of this property, if we do not insist on $\sigma$-additivity.\footnote{There is a controversy among philosophers about the type of additivity that should be required for a probability measure. We take the side of finite additivity in this debate, much in line with the critique of $\sigma$-additivity by \citet{definetti}. Notable arguments for $\sigma$-additivity in this debate are found in the work of, for example, \citet{skyrms} and \citet{williamson}. }

In addition to the above, there is also the following argument that uses the supertask itself as an argument against $\sigma$-additivity, analogue to the argument of \citet{laraudogoita}. Any probability function consistent with the supertask, has to assign probability zero to the event that the final ball is ball $k$. Any $\sigma$-additive probability function, however, cannot assign zero probability to this event for every $k$: because there are countably many balls and since the outcome space is the disjoint union of all these events, it then also should have probability zero, which would contradict that the outcome space has probability $1$ by definition. So if we want that our model to allow at least one probability function to be consistent with the supertask, we cannot insist on $\sigma$-additivity.

Within our model, we prove that for any $p \in [0,1]$ and any set $A \subseteq \mathbb{N}$ such that $A$ and $A^c$ are infinite, there is a probability function satisfying the constraint, assigning probability $p$ to the final ball being in $A$. To show the idea behind the proof, we fix the outcome
\begin{equation}
\{1\},\{1,2\},\{1,2,3\},\{1,2,3,4\},\ldots \;.
\end{equation}
Consider the probability function that assigns probability $1/2$ to the two outcomes that are identical to our fixed outcome up to god $1$, i.e.
\begin{equation}
\begin{aligned}
\{1\},\{1,2\},\{1,2,3\},\{1,2,3,4\},\ldots & \;\mathrm{and}\\
\{2\},\{1,2\},\{1,2,3\},\{1,2,3,4\},\ldots & \;. \\
\end{aligned}
\end{equation}
Since what happens up to god $1$ is deterministic ($\{1,2\},\{1,2,3\},...$ with probability $1$) under this probability function, the constraint for $k=1$ entails that god $1$ removes ball $1$ or $2$ both with probility $1/2$, which is exactly what happens. So this probability function satisfies the constraint for $k=1$, but not for all other $k$. Now we consider the probability function that assigns probability $1/6$ to the outcomes that are identical to our fixed outcome up to god $2$, i.e.
\begin{equation}
\begin{aligned}
\{1\},\{1,2\},\{1,2,3\},\{1,2,3,4\},\ldots & \;, \\
\{2\},\{1,2\},\{1,2,3\},\{1,2,3,4\},\ldots & \;,\\
\{2\},\{2,3\},\{1,2,3\},\{1,2,3,4\},\ldots & \;,\\
\{3\},\{2,3\},\{1,2,3\},\{1,2,3,4\},\ldots & \;,\\
\{1\},\{1,3\},\{1,2,3\},\{1,2,3,4\},\ldots & \;\mathrm{and}\\
\{3\},\{1,3\},\{1,2,3\},\{1,2,3,4\},\ldots & \;.\\
\end{aligned}
\end{equation}  
What happens up to god $2$ is now deterministic and since every possibility for the actions of gods $1$ and $2$ have the same probability, the constraint is satisfied for $k=1$ and $k=2$. By continuing this process, we get probability functions that satisfy the constraint for more and more $k$. In this way, in the limit, we get a probability function that satisfies the constraint for every $k$ (see Section \ref{sec:model} for the precise construction of this limit probability function). If what happens up to god $k$ is deterministic and everything after that happens with the same probability, then the probability of the final ball being even is given by the fraction of even balls after god $k+1$ has removed a ball. Since our fixed outcome is  $\{1\},\{1,2\},\{1,2,3\},\ldots$, the probability of the final ball being even converges to $1/2$ as $k$ grows to infinity. Hence the probability that the final ball is even is $1/2$.

The crux is that we can start with \emph{any} fixed outcome. Suppose we want a probability function that assigns probability $1/3$ to the final ball being even. Then we construct an appropriate fixed outcome in the following way. We start with $\{1\}$. Because the fraction of even numbers is now $0$, we add an even number for the next element of the sequence: $\{1\},\{1,2\}$. The fraction of even numbers is now $1/2$, which is greater than $1/3$, so we add an odd number:  $\{1\},\{1,2\},\{1,2,3\}$. Now the fraction of even numbers is exactly $1/3$, so we can add any number: $\{1\},\{1,2\},\{1,2,3\},\{1,2,3,4\}$. Because the fraction is now again greater than $1/3$, we add an odd number and so on. So if we fix the outcome
\begin{equation}
\begin{aligned}
& \{1\},\{1,2\},\{1,2,3\},\{1,2,3,4\},\{1,2,3,4,5\},\{1,2,3,4,5,7\}, \\
& \{1,2,3,4,5,6,7\},\{1,2,3,4,5,6,7,9\},\{1,2,3,4,5,6,7,9,11\},\ldots \;,
\end{aligned}
\end{equation}
then the probability of the final ball being even converges to $1/3$ instead of $1/2$. The same works for any given $p \in [0,1]$ other than $1/3$. Further, the set of even numbers is not special here. The only thing we use to construct the appropriate outcome is that at any time we can add an even or an odd number, in other words: both the set of even numbers and its complement have to be infinite. Hence the same works for any infinite set of natural numbers with an infinite complement, as desired.

To see what happens for finite set of natural numbers, we look at the probability of the final ball being a specific ball, for example ball $7$. Suppose that ball $7$ is still in the urn when god $k$ takes the urn. For ball $7$ to make it to the end, god $k$ to $1$ have to remove some other ball than $7$. They do that with respectively probability $k/(k+1)$ to $1/2$, making the probability that $7$ is the remaining ball
\begin{equation}
\frac{k}{k+1} \cdot \frac{k-1}{k} \cdot \ldots \cdot \frac{2}{3} \cdot \frac{1}{2} = \frac{1}{k+1}.
\end{equation} 
So the probability ball $7$ survives from god $k$ to $1$, is $1/(k+1)$. For ball $7$ to be the remaining ball, however, it has to survive \emph{all} gods, which implies that the probability of that happening has to be smaller than $1/(k+1)$ for every $k$. Hence the probability the final ball is ball $7$ can only be zero. In the same way, any specific ball has probability zero of being the final ball. And since our probability function is finitely additive, the probability of the final ball being from any finite set of balls, is zero. 

\section{The mathematical model}
\label{sec:model}

In this section, we make everything mathematically precise. To do that, we start with some definitions. Write for $k \in \mathbb{N}$
\begin{equation}
\begin{aligned}
\mathbb{N}_k & := \{ A \subseteq \mathbb{N} \;:\; |A|=k\} \\
\Omega_k & := \{ (B_k,B_{k+1},\ldots) \;:\; B_j \in \mathbb{N}_j, \; B_j \subset B_{j+1} \;\mathrm{and}\; \cup_j B_j = \mathbb{N}  \}.
\end{aligned}
\end{equation}
Elements of $\Omega_k$ represent what happens up until god $k$. The space $\Omega:=\Omega_1$ is our outcome space. On this space, we define the following two functions. For $B = (B_{1},B_{2},...) \in \Omega$ with $B_1=\{n\}$, we define 
\begin{equation}
\begin{aligned}
H_k(B) & := (B_k,B_{k+1},\ldots) \in \Omega_{k+1} \\
R(B) & := n.
\end{aligned}
\end{equation}
The function $H_k$ maps an outcome to what happened up until god $k$ and $R$ gives us the remaining ball at the end of the process. 

For $B = (B_{k+1},B_{k+2},...) \in \Omega_{k+1}$, $S \subseteq \Omega_k$ and $1 \leq j \leq k+1$ we set
\begin{equation}
\begin{aligned}
N(S; B)) & := |\{  b \in B_{k+1} \;:\; (B_{k+1} \setminus \{b\},B_{k+1},B_{k+2},\ldots) \in S \}| \\
P_j(S) & := \left\{ B \in \Omega_{k+1} \;:\; N(S; B)=j \right\} \subseteq \Omega_{k+1}
\end{aligned}
\end{equation}
In words, the set $P_j(S)$ gives us everything that can happen up until god $k+1$ such that there are precisely $j$ of the $k+1$ balls that upon removal give us an element of $S$. 

We consider the  measure space $(\Omega, \mathcal{P}(\Omega))$. On this measure space we consider probability functions as candidates. The constraint on probability function $\mu$ that we obtain from the description of the supertask is that for every $k$, $S \subseteq \Omega_k$ and $T \subseteq P_j(S)$ such that $\mu(H_{k+1} \in T)>0$ we have
\begin{equation}
\label{eq:constraint}
\mu(H_k \in S | H_{k+1} \in T) = \frac{j}{k+1}.
\end{equation}
We note that using (\ref{eq:constraint}) we can determine $\mu(H_k \in S | H_{k+1} \in T)$ for every $k$, $S \subseteq \Omega_k$ and $T \subseteq \Omega_{k+1}$ such that $\mu(H_{k+1} \in T)>0$, by writing
\begin{equation}
\begin{aligned}
\mu(H_k \in S | H_{k+1} \in T) & = \sum_{j=0}^{k+1} \frac{\mu(H_k \in S; H_{k+1} \in T \cap P_j(S))}{\mu(H_{k+1} \in T)} \\
&  = \sum_{j=1}^{k+1} \frac{j}{k+1} \frac{\mu(H_{k+1} \in T \cap P_j(S))}{\mu(H_{k+1} \in T)}.
\end{aligned}
\end{equation}
We write $\mathcal{M}$ for the set of probability functions that satisfy (\ref{eq:constraint}) and prove the following theorem.

\begin{theorem}
\label{thm:theoremA}
Let $A \subseteq \mathbb{N}$. Then
\begin{equation}
\{ \mu(R \in A) \;:\; \mu \; \in \mathcal{M} \} = \left\{
	\begin{array}{lll}
		\{0\}  & \;\; \mbox{if } A\mathrm{\;is\;finite} \\
		\{1\} & \;\;\mbox{if } A^c\mathrm{\;is\;finite} \\
		{[0,1]} & \;\;\mathrm{otherwise}
	\end{array}
\right..
\end{equation}
\end{theorem}

We present the proof of Theorem \ref{thm:theoremA} in three steps.

\textbf{Step 1} We construct for every $Z \in \Omega$, a probability function $\alpha_Z \in \mathcal{M}$ based on $Z$.

To construct a probability function such that (\ref{eq:constraint}) is satisfied, we fix a $Z:= (Z_1,Z_2,Z_3,\ldots) \in \Omega$ and set
\begin{equation}
F_n := \{ H_n = (Z_n,Z_{n+1},\ldots) \} \subseteq \Omega.
\end{equation}
For $S \subseteq \Omega = \Omega_1$ we define the sequence $x(S) \in [0,1]^\infty$ by
\begin{equation}
x_n(S) := \frac{|S \cap F_n|}{|F_n|} = \frac{|S \cap F_n|}{n!} ,
\end{equation}
which gives the density of $S$ with respect to $F_n$. Now, we let $n$ go to infinity. To do that, we extend the limit operator on all convergent sequences to a linear operator $L$ on all bounded sequences. Such an extension can be constructed by using a free ultrafilter on $\mathbb{N}$. The existence of such a free ultrafilter is guaranteed by the Boolean Prime Ideal Theorem, which can not be proven in ZF set theory, but is weaker than the axiom of choice \citep{halpern}.\footnote{The existence of a atomfree or nonprincipal (i.e. every singleton has measure zero) probability function defined on the power set of $\mathbb{N}$ cannot be established in ZF alone \citep{solovay}. Consequently, a version of the axiom of choice is always necessary to construct an element of $\mathcal{M}$.} We define $\alpha_Z: \mathcal{P}(\Omega) \rightarrow [0,1]$ by 
\begin{equation}
\alpha_Z(S) := L(x(S)).
\end{equation}

Now we show that $\alpha_Z \in \mathcal{M}$. For $S \subseteq \Omega_k$, $T \subseteq P_j(S)$ and $k<n$, we have
\begin{equation}
\label{eq:sumtotal}
\begin{aligned}
 j |\{ H_{k+1} \in T\} \cap F_n| & = \sum_{B \in T} N(S,B) |\{H_{k+1}=B\} \cap F_n |  \\
& = \sum_{B' \in S} \sum_{B \in T} N(\{B'\},B) |\{H_{k+1}=B\} \cap F_n|  \\
& = \sum_{B' \in S} \sum_{B \in T} (k+1) |\{H_k=B'; H_{k+1}=B\} \cap F_n| \\
& = (k+1) |\{H_k \in S; H_{k+1} \in T\} \cap F_n|.
\end{aligned}
\end{equation}
Dividing by $n!(k+1)$ on both sides in (\ref{eq:sumtotal}) we find that
\begin{equation}
x_n(\{H_k \in S; H_{k+1} \in T\}) = \frac{j}{k+1} x_n(\{ H_{k+1} \in T \}).
\end{equation}
Since $L$ is linear, we get that
\begin{equation}
\frac{\alpha_Z\left( H_k \in S; H_{k+1} \in T \right)}{\alpha_Z(H_{k+1} \in T)} =  \frac{j}{k+1} 
\end{equation}
for every $k$, $S \subseteq \Omega_k$ and $T \subseteq P_j(S)$ with $\alpha_Z(H_{k+1} \in T)>0$, which is what we wanted to show.

\textbf{Step 2} We show that for every infinite $A \subseteq \mathbb{N}$ such that $A^c$ is also infinite and every $p \in [0,1]$ there is a $Z \in \Omega$ such that $\alpha_Z(R \in A)=p$. 

Let $A \subseteq \mathbb{N}$ be such that both $A$ and $A^c$ are infinite. Let $p \in [0,1]$ be given. If we choose $Z \in \Omega$ such that
\begin{equation}
\label{eq:densitytop}
\frac{|Z_k \cap A|}{k} \rightarrow p
\end{equation}
as $k \rightarrow \infty$, then $\alpha_Z(R \in A) = p$. We construct $Z$ by recursion, using as basic rule that we add an element from $A$ if the density is too high and add an element from $A^c$ if the density is too low. To make sure that we add every natural number somewehere in the process and thus $Z \in \Omega$, we make one exception to this rule: for every $j$, we add an element from $A$ in the $j^2$-th step and an element from $A^c$ in the $j^2+1$-th step. This exception, however, does not influence the limit density we obtain since the fraction of square numbers in $[0,k]$ converges to zero as $k \rightarrow \infty$. We set
\begin{equation}
\begin{aligned}
Z_1 & := \{1\} \\
a_k & := \min A \cap Z_k^c \\
b_k & := \min A^c \cap Z_k^c \\
Z_{k+1} & :=  \left\{
	\begin{array}{lll}
		Z_k \cup \{a_k\} & \;\; \mbox{if } k \not\in \{j^2,j^2+1\} \mathrm{\;for\;some\;}j \mathrm{\;and\;} \frac{|Z_k \cap A|}{k} \leq p \\
		Z_k \cup \{b_k\} & \;\; \mbox{if } k \not\in \{j^2,j^2+1\} \mathrm{\;for\;some\;}j \mathrm{\;and\;} \frac{|Z_k \cap A|}{k} > p \\
		Z_k \cup \{a_k\} & \;\;\mbox{if } k=j^2 \mathrm{\;for\;some\;}j \\
		Z_k \cup \{b_k\} & \;\;\mbox{if } k=j^2+1 \mathrm{\;for\;some\;} j \\
	\end{array}
\right..
\end{aligned}
\end{equation}
Notice that $a_k$ and $b_k$ are well defined for every $k$ because $A$ and $A^c$ are both infinite. Also note that $Z=(Z_1,Z_2,Z_3,...) \in \Omega$ and (\ref{eq:densitytop}) is satisfied, as desired.

\textbf{Step 3} We show that for all $\mu \in \mathcal{M}$ and finite $A \subseteq \mathbb{N}$ we have $\mu(R \in A)=0$.

Let $a \in \mathbb{N}$ and $\mu \in \mathcal{M} \not= \emptyset$. If we take
\begin{equation}
S_k := \{ (B_k,B_{k+1},\ldots) \in \Omega_k \;:\; a \in B_k \},
\end{equation}
we get from (\ref{eq:constraint}) that
\begin{equation}
\label{eq:removalofa}
\mu(H_k \in S_k ) = \frac{k}{k+1} \mu(H_{k+1} \in S_{k+1}).
\end{equation}
This implies that $\mu(R=a) = \mu(H_1 \in S_1) = \frac{1}{k} \mu(H_k \in S_k)$ for every $k$, so $\mu(R=a) \leq \frac{1}{k}$ for every $k$. Hence $\mu(R=a)=0$. Since $\mu$ is finitely additive, we have $\mu(R \in A)=0$ for every finite $A \subseteq \mathbb{N}$.
\qed

\section{Discussion}
\label{sec:discussion}

We have given a mathematical model for the supertask of Hansen. Our conclusions depend of course on the choice for this model: another model could give different results. We think, however, that the present model is not controversial. From our model, we have concluded that there are many probability functions that are consistent with the supertask. We have expressed this indeterminacy in terms of the different probabilities the probability functions assign to the final ball being in some subset $A \subseteq \mathbb{N}$. If $A$ is infinite and its complement is also infinite, then any probability $p$ is assigned to the final ball being in $A$ by some probability function consistent with the supertask.

This conclusion shows that the description of the supertask as presented, just does not give us enough information to say \emph{anything} about the probability of the final ball being in a set of the type described above. As a consequence, we should be very careful about this indeterminacy when thinking about the supertask. At one point, Hansen mentions the `emperical distribution' obtained by the gods if they keep repeating the supertask. Assuming there is the same underlying probability function every time they perform the supertask, of course one can, at least in theory, estimate this probability function through performing experiments. However, the description of the supertask does not give enough information to decide what that underlying probability function is.

One would hesitate to call a probability function assigning probability $1$ to the even numbers and probability $0$ to the odd numbers, a `uniform' probability function. A reasonable notion of uniformity at least includes that every residue class modulo $m$ has probability $1/m$ (giving for $m=2$ that the even numbers have probability $1/2$).\footnote{This is the weakest notion of uniformity considered by \citet{schirokauerkadane}. } Our model shows, however, that the final ball being even can get any probability. This means that the notion of uniformity emerging from the supertask is weaker than any reasonable notion of uniformity.

Our analysis also shows that for any  $n \in \mathbb{N}$, a probability function consistent with the supertask necessarily has to assign probability zero to the final ball being $n$. Altough this property of being atomfree is not a reasonable notion of `uniformity', it is some form of `fairness' nonetheless. This form of fairness is sufficient for the point Hansen wants to make about the supertask, namely that given any possible final ball $k$, the probability is $1$ that the final ball is bigger than $k$. In other words: the final ball seems always unexpectedly low. This is certainly not an inconsistency, but Hansen does call it an `absurdity'.

We do not want to present an argument here about whether this is indeed an absurdity or not, but only to add clarity to that discussion by a proper analysis of the supertask. We do, however, want to point out that the alleged absurd property of the probability function on the final ball is not unique. Also it is not restricted to probability functions. A typical $\sigma$-additive probability measure with this property is the following. Let $\omega_1$ be the first uncountable ordinal number and let $\mathcal{F}$ be the $\sigma$-algebra of subsets $F \subseteq \omega_1$ such that either $F$ or $F^c$ is countable. Let $\mu : \mathcal{F} \rightarrow [0,1]$ be given by
\begin{equation}
\mu(F) := \left\{
	\begin{array}{lll}
		0 & \;\; \mbox{if } F \mathrm{\;is\;countable} \\
		1 & \;\; \mbox{if } F \mathrm{\;is\;uncountable}
	\end{array}
\right..
\end{equation}
Then $\mu$ is a $\sigma$-additive probability measure. If $\mu$ models a lottery on $\omega_1$, then completely analogous to the situation on $\mathbb{N}$, for any outcome $\alpha \in \omega_1$ we have $\mu(\{ \beta \in \omega_1 \;:\; \beta > \alpha\})=1$. Accepting the Axiom of Choice, this property also translates to picking a random number from $[0,1]$. By the well ordering theorem, there exists a well ordering $\preceq$ of $[0,1]$. If we model the experiment with  Lebesgue measure $\lambda$ on $[0,1]$, for any outcome $x \in [0,1]$, we again have $\lambda(\{y \in [0,1]\;:\; y \succ x\})=1$.

\subsection*{Acknowledgement}
I would like to thank my supervisor Ronald Meester for all the very helpful discussions.






\begin{flushright}
Timber Kerkvliet\\
VU University Amsterdam\\
\href{mailto:t.kerkvliet@vu.nl}{t.kerkvliet@vu.nl}\\

\end{flushright}

\end{document}